\documentclass[10pt]{amsart}
\title[Multi-marginal optimal transport with quadratic costs and uniform discrete marginals]{\textbf{On the existence of Monge solutions to multi-marginal optimal transport with quadratic cost and uniform discrete marginals
}}
\date{\today}
\author{Pedram Emami}
\email{emami1@ualberta.ca, pass@ualberta.ca}
\author{Brendan Pass}\thanks{BP is pleased to acknowledge support from Natural Sciences and Engineering
Research Council of Canada Grant 04658-2018. The work of PE was completed in partial fulfillment of the
requirements for a doctoral degree in mathematics at the University of Alberta.}
\address{University of Alberta, Edmonton, Alberta, Canada}

\usepackage{amsmath,amsthm,amssymb,mathrsfs,bm}
\usepackage{float,caption}
\usepackage{geometry,graphicx,tabularx,booktabs}
\usepackage{xcolor,colordvi,color}
\usepackage{hyperref}

	\definecolor{heavyred}{cmyk}{0,1,1,0.25}
	\definecolor{heavyblue}{cmyk}{1,1,0,0.25}
	\hypersetup{
	  pdftitle=,
	  pdfpagemode=UseOutlines,
	  pageanchor=false,
	  citebordercolor=0 0 1,
	  colorlinks=true,
	  allcolors=heavyred,
	  breaklinks=true,
	  pdfauthor=,
	  pdfpagetransition=Dissolve,
	}
	\newtheoremstyle{thmStyle}{}{}{\it}{}{\bfseries}{.}{ }{}
		\newtheoremstyle{corStyle}{}{}{\it}{}{\bfseries}{.}{ }{}
	\newtheoremstyle{prfStyle}{}{}{}{}{\it\bfseries}{.}{ }{}
	\newtheoremstyle{exStyle}{}{}{}{}{\bfseries}{.}{ }{}	
	{
	  \theoremstyle{thmStyle}
	  \newtheorem{theorem}{Theorem}
	}
	{
	\theoremstyle{corStyle}
	\newtheorem{corollary}{Corollary}
}
	{
   	  \theoremstyle{prfStyle}
	  \newtheorem*{prf}{Proof}
	}
	{
   	  \theoremstyle{exStyle}
	  \newtheorem{example}{Example}
	}

\def\R{\mathbb{R}}
\def\(#1\){\left(#1\right)}
\def\]#1\]{\left\{#1\right\}}
\def\|#1\|{\left|#1\right|}
\def\half{\frac{1}{2}}
\def\trd{\frac{1}{3}}

\def\etal{\@ifnextchar,{\textit{et al.}}{\@ifnextchar.{\textit{et al}}{\textit{et al.\@ }}}}

\newcommand*{\mapSubTwo}[1][]{\def\mapOptTwo{#1}%
\ifx\mapOptOne\empty%
	\def\mapOptOne{X}%
\fi
\ifx\mapOptTwo\empty%
	\def\mapOptTwo{Y}%
\fi
\mapMain:{\mapOptOne}\rightarrow\mapOptTwo%
}

\begin{document}
\begin{abstract}
A natural and important question in multi-marginal optimal transport is whether the \emph{Monge ansatz}
is justified; does there exist a solution of Monge, or deterministic, form?  We address this question for
the quadratic cost when each marginal measure is $m$-empirical (that is, uniformly supported on $m$
points). By direct computation, we provide an example showing that the ansatz \emph{can fail} when the
underlying dimension $d$ is $2$, number of marginals $N$ to be matched is $3$ and the size $m$ of their
supports is $3$. As a consequence, the set of $m$-empirical measures is not barycentrically convex when
$N \geq 3$, $d \geq 2$ and $m \geq3$.  It is a well known consequence of the Birkhoff-von Neumann Theorem
that the Monge ansatz holds for $N=2$, standard techniques show it holds when $d=1$, and we provide a
simple proof here that \emph{it holds whenever $m=2$}.  Therefore, the $N$, $d$ and $m$ in our
counterexample are as small as possible.
\end{abstract}
\maketitle

\section{Introduction}
Given probability measures  $\mu_1,\mu_2,...,\mu_N$ (known as the \emph{marginals}) on respective bounded
domains $X_1,X_2,...,X_N \subseteq \mathbb{R}^d$ and a cost function $c: X_1 \times X_2\times...\times
X_N \rightarrow \mathbb{R}$, the multi-marginal optimal transport problem is to minimize
\begin{equation}\label{eqn: OT problem}
\inf_{\gamma\in\Pi(\mu_1,\mu_2,...,\mu_N)}\int_{X_1\times X_2 \times ...\times
X_N}c(x_1,x_2,...,x_N)\,d\gamma,
\end{equation}
over the set $\Pi(\mu_1,\mu_2,...,\mu_N)$ of probability measures on $X_1\times X_2 \times ...\times X_N$
whose marginals are the $\mu_i$. In particular, the $N=2$ case is the classical optimal transport
problem, which represents a very active area of modern mathematics, with many applications, reviewed in,
for example \cite{villani2021topics, santambrogio2015optimal,Villani2009}.  The more general case, $N\geq
3$, has recently emerged as an important problem in its own right, owing largely to it own wide variety
of applications; see, for instance, the survey papers \cite{Pass2015,DiMarinoGerolinNenna2017}.

A central question in this theory is whether minimizers concentrate on the graphs of functions over the
first variable $x_1$; solutions with this structure are known as \emph{Monge solutions}. Restricting the
minimization to measures arising from such mappings, or, adopting the terminology of
\cite{Friesecke2019}, making the \emph{Monge ansatz}, when justified, amounts to a spectacular
dimensional reduction and therefore greatly simplifies the problem computationally. In the classical
$N=2$ case, it is well known that if $\mu_1$ is absolutely continuous with respect to Lebesgue measure,
then, for cost satisfying mild structural conditions, the solution to \eqref{eqn: OT problem} is unique
and of Monge type \cite{Caffarelli1996,GangboMcCann1996,Levin1999}.  Central among such costs is the
\emph{quadratic cost}, $c(x_1,x_2) = |x_1-x_2|^2$, which induces a metric, known as the \emph{Wasserstein
distance}, on the set of probability measures, defined by
$W_2(\mu_1,\mu_2):=\sqrt{\inf_{\gamma\in\Pi(\mu_1,\mu_2)}\int_{X_1\times X_2}|x_1-x_2|^2\,d\gamma}$;
Monge solution and uniqueness results for this cost were first proven by Brenier \cite{Brenier1987}. The
situation is much more subtle and less well understood when $N \geq 3$.  However, by now many examples of
cost functions for which unique Monge solutions exist for absolutely continuous marginals
\cite{Carlier2003,Heinich2002,ColomboDePascaleDiMarino2015,PassVargasJimenez2021,PassVargas-Jimenez2023},
as well as many for which they do not \cite{CarlierNazaret2008,Pass2013,Pass2012a}, have been discovered,
and some general (though very strong) conditions on $c$ ensuring this structure have been developed
\cite{Pass2012,KimPass2014,PassVargasJimenez2023b}.  Of special note is the natural, multi-marginal
extension of the quadratic cost,
\begin{equation}\label{eqn: multi-marginal quadratic cost} 
c(x_1,x_2,...,x_N) = \sum_{i,j=1}^N|x_i-x_j|^2,
\end{equation}
or equivalently, $c(x_1,x_2,...,x_N) = -\|\sum_{i=1}^Nx_i\|^2$, for which existence of Monge solutions
and uniqueness were proven by Gangbo and Swiech \cite{GangboSwiech1998}.  Solving \eqref{eqn: OT
problem} with this cost is known to be equivalent to finding a \emph{Wasserstein barycenter} (with equal
weights on the measures), which, as introduced by Agueh-Carlier \cite{agueh2011barycenters}, is a
minimizer among probability measures $\nu$ on $\mathbb{R}^d$ of
\begin{equation}\label{eqn: Wasserstein barycenter}
\nu\mapsto\sum_{i=1}^N W_2^2(\nu,\mu_i).
\end{equation}
Wasserstein barycenters are Fr\'{e}chet means on the Wasserstein space, and are an extremely active area
of current research, due to many applications in statistics, image processing and data science, among
other areas (see, for instance, \cite{Rabinetal2012,Hoetal2017,YangTabak2022,StaibClaiciSolomonJegelka17}
and the survey in the monograph \cite{PeyreCuturi2019}).

In this short note we are interested in discrete, rather than absolutely continuous, marginals $\mu_i$.
In general, optimal couplings between discrete measures need not be unique or of Monge type; indeed, in
general, there may not even exist a single coupling $\gamma \in \Pi(\mu_1,\mu_2,...,\mu_N)$ which is of
Monge type. However, such \emph{transport maps} at least exist when all the marginals are uniformly
distributed on a common number $m$ of points, as is often the case in applications. We will call such
measures \emph{$m$-empirical measures}.  In this case, when $N=2$, the \emph{existence} of Monge
solutions for \emph{any} cost function $c$ is an easy consequence of the Birkhoff-von Neumann Theorem,
although solutions to \eqref{eqn: OT problem} are not unique in general. An immediate consequence is that
the set of $m$-empirical measures is \emph{geodesically convex} in Wasserstein space; that is, any two
elements can be connected by a Wasserstein geodesic (often called a displacement interpolant) which lies
in this set.

On the other hand, even for $m$-empirical measures, Monge solutions for $N \geq 3$ do not exist for some
cost functions. Indeed, as exposed in \cite{Friesecke2019}, when $N \geq 3$, the convex polytope
$\Pi(\mu_1,...,\mu_N)$ has extreme points which are not of Monge type, unlike in the $N=2$ case.  This
immediately implies the existence of cost functions for which Monge solutions do not exist, and
\cite{Friesecke2019} provides examples with relevance in physics for which this occurs.  Intuition
provided there, however, suggests that it is the \emph{repulsive nature of certain costs} which leads to
non-Monge optimizers for $m$-empirical measures, in agreement with known results in the continuous case.
It remains reasonable to expect existence of Monge solutions for the quadratic cost function
\eqref{eqn: multi-marginal quadratic cost}. It is similarly natural to conjecture that $m$-empirical
marginals always have $m$-empirical Wasserstein barycenters (that is, that the set of $m$-empirical
measures is \emph{barycentrically convex}); we note that upper bounds on the size of the support of
barycenters are known for general discrete measures \cite{anderes2016discrete}, but the bound obtained
there is $N(m-1)+1$ in the case of $m$-empirical marginals, much larger than $m$.

We exhibit here, by direct computation, an example showing that Monge solutions may not exist for $N\geq
3$ even for the quadratic cost. An immediate consequence is that the set of $m$-empirical measures is not
barycentrically convex; there exist $m$-empirical marginals $\mu_i$ for which there does not exist a
Wasserstein barycenter within the same set.

Our counterexample is for three measures ($N=3$) in two dimensions ($d=2$) supported on three points
($m=3$).  These are the smallest values of these three parameters for which  such counterexamples can
exist.  Indeed, for $N=2$, Monge solutions can always be found as discussed above.  When $d=1$, optimal
couplings are well known to be monotone increasing, which implies they are of Monge form and unique.
Finally, we present here a simple proof that there always exists a Monge solution for any collection of
$2$-empirical marginals, for \emph{any} values of $d$ and $N$.

\section{A counterexample to Monge solutions for marginals with three-point support in $\R^2$} 
We define the set of $m$-empirical measures on $\mathbb{R}^d$ as
$E_m:=\{\frac{1}{m}\sum_{i=1}^m\delta_{x_i}: x_1,..,x_m \in \mathbb{R}^d\}$.

Here we present a counterexample to the existence of Monge solutions for $3$ marginals supported
uniformly on $3$ points in $\mathbb{R}^2$.
\begin{example}\label{ex: counterexample}
Let $\mu_1,\mu_2,\mu_3 \in E_3$ be the $3$-empirical measures supported on $3$ points in $\R^2$ defined by
\begin{eqnarray}
\mu_1&=\trd\(\delta_{(0.4417,-4.7665)}+\delta_{(-0.27748,1.0397)}+\delta_{(1.4826,4.7896)}\),\nonumber\\
\mu_2&=\trd\(\delta_{(-2.1054,-3.9784)}+\delta_{(3.5763,-1.8988)}+\delta_{(3.328,-1.558)}\),\nonumber\\
\mu_3&=\trd\(\delta_{(-3.6728,0.23451)}+\delta_{(1.6988,-2.2917)}+\delta_{(-1.1644,-2.386)}\).\label{eqn:
marginals for counterexample}
\end{eqnarray} 
The optimal transport problem \eqref{eqn: OT problem} with cost \eqref{eqn: multi-marginal quadratic
cost} can easily be solved directly by linear programming, yielding the optimal measure
$\gamma=\frac{1}{6}[\delta_{x^{113}} +
\delta_{x^{122}}+\delta_{x^{211}}+\delta_{x^{223}}+\delta_{x^{331}}+\delta_{x^{332}}]$, with an optimal
cost of $\mathbf{68.027}$.  Here, $x^{\alpha^1\alpha^2\alpha^3} := (x_1^{\alpha^1},
x_2^{\alpha^3},x_3^{\alpha^3}) \in (\mathbb{R}^2)^3$, where the elements $x_i^{\alpha^i}\in\mathbb{R}^2$
in the support of each $\mu_i$ are ordered as in \eqref{eqn: marginals for counterexample}: for example,
$x^{113} =((0.4417,-4.7665),(-2.1054,-3.9784),(-1.1644,-2.386))$.

On the other hand, potential Monge solutions arise from transport maps, which correspond to measures of
the form $\gamma=\frac{1}{3}[\delta_{x^{1,\sigma_2(1),\sigma_3(1)}} +
\delta_{x^{2,\sigma_2(2),\sigma_3(2)}}+\delta_{x^{3,\sigma_2(3),\sigma_3(3)}}]$ for permutations
$\sigma_2$ and $\sigma_3$ of the set $\{1,2,3\}$ of indices.  One can easily check that the
\textit{minimal Monge cost} (MMC), obtained by minimizing the transport cost among the $(N!)^{m-1}
=(3!)^2 =36$ transport maps, is obtained by $\gamma=\frac{1}{3}[\delta_{x^{113}} +
\delta_{x^{222}}+\delta_{x^{331}}]$, yielding a strictly higher total cost of $\mathbf{68.065}$.  It
follows that the problem \eqref{eqn: OT problem} with quadratic cost \eqref{eqn: multi-marginal quadratic
cost} and marginals \eqref{eqn: marginals for counterexample} does not admit a Monge solution.

Since a $3$-empirical Wasserstein barycenter of $\mu_1,\mu_2$ and $\mu_3$ would immediately yield a Monge
solution to \eqref{eqn: OT problem} via the equivalence between the problems identified in
\cite{agueh2011barycenters}, this also implies that there cannot exist a $3$-empirical Wasserstein
barycenter $\nu \in E_3$ of these marginals.
\end{example}
\section{Existence of Monge solutions for uniform empirical measures with two-point supports in $\R^d$}
We now present a simple proof that Monge solutions exist for $2$-empirical marginals.

\begin{theorem}
Let $\mu_i=\half\(\delta_{x_i^1}+\delta_{x_i^2}\) \in E_2$ be $2$-empirical measures on $\mathbb{R}^d$
for $i=1,...,N$.  Then problem \eqref{eqn: OT problem} with cost \eqref{eqn: multi-marginal quadratic
cost} has a Monge solution.
\end{theorem}
\begin{prf}
Since the transform $x_i \rightarrow x_i-\bar x_i$, where each $x_i :=(x_i^1+x_i^2)/2$ is the mean of
$\mu_i$, shifts the values of the functional \eqref{eqn: OT problem} but does not change its minimizer,
we can assume without loss of generality that the mean of each $\mu_i$ is $0$; that is, $x_i^2=-x_i^1$.
Now, choose $\alpha_i \in \{1,2\}$ for $i=1,2...,N$, to maximize
\begin{equation*}
M = \max_{\alpha_i \in \{1,2\}: i=1,2,...,N}\|\sum_{i=1}^nx_i^{\alpha_i}\|^2.
\end{equation*}

Then set $\hat x =(x_1^{\alpha^1},x_2^{\alpha^2},...,x_N^{\alpha^N}) \in (\mathbb{R}^{d})^N$ and $\gamma
= \frac{1}{2}[\delta_{\hat x} +\delta_{-\hat x}]$.  It is clear that $\gamma \in \Pi(\mu_1,...,\mu_N)$
and that $\gamma$ is supported on the graph of the transport map $x_1^{\alpha^1} \mapsto
(x_2^{\alpha^2},...,x_N^{\alpha^N})$, $-x_1^{\alpha^1}\mapsto (-x_2^{\alpha^2},...,-x_N^{\alpha^N})$.
Since $\|\sum_{i=1}^nx_i\|^2=M$ $\gamma$ a.e, we have
\begin{equation*}
-\int_{(\mathbb{R}^{d})^N} \|\sum_{i=1}^nx_i\|^2 d\gamma =-M,
\end{equation*}
whereas, for any other $\tilde \gamma \in \Pi(\mu_1,...\mu_N)$ we have
\begin{equation*}
-\int_{(\mathbb{R}^{d})^N} \|\sum_{i=1}^nx_i\|^2 d\tilde \gamma \geq -M,
\end{equation*}
establishing the optimality of the Monge coupling $\gamma$ for the cost function
$c(x_1,x_2,...,x_N)=-\|\sum_{i=1}^nx_i\|^2$, which is well known to be equivalent to \eqref{eqn:
multi-marginal quadratic cost} \cite{GangboSwiech1998}.
\end{prf}
% ------------------
\begin{corollary}
Let $\mu_i=\half\(\delta_{x_i^1}+\delta_{x_i^2}\) \in E_2$ be $2$-empirical measures on $\mathbb{R}^d$
for $i=1,...,N$.  Then there exists a  $2$-empirical Wasserstein barycenter $\nu \in E_2$ of the $\mu_i$.
\end{corollary}
\begin{prf}
Let $x^{\pm} :=\pm\frac{1}{N}\sum_{i=1}^Nx_i^{\alpha_i} \in \mathbb{R}^2$, where the $\alpha_{i} \in
{1,2}$ are as in the previous proof.  The equivalence between multi-marginal optimal transport and the
Wasserstein barycenter established in Proposition 4.1 of \cite{agueh2011barycenters} then implies that
$\frac{1}{2}[\delta_{x^+} +\delta_{x^-}] \in E_2$ is a Wasserstein barycenter.
\end{prf}
\section{Discussion}
\begin{figure}[H]
\centering
\captionsetup{font=small,justification=centering}
\includegraphics[width=0.75\textwidth]{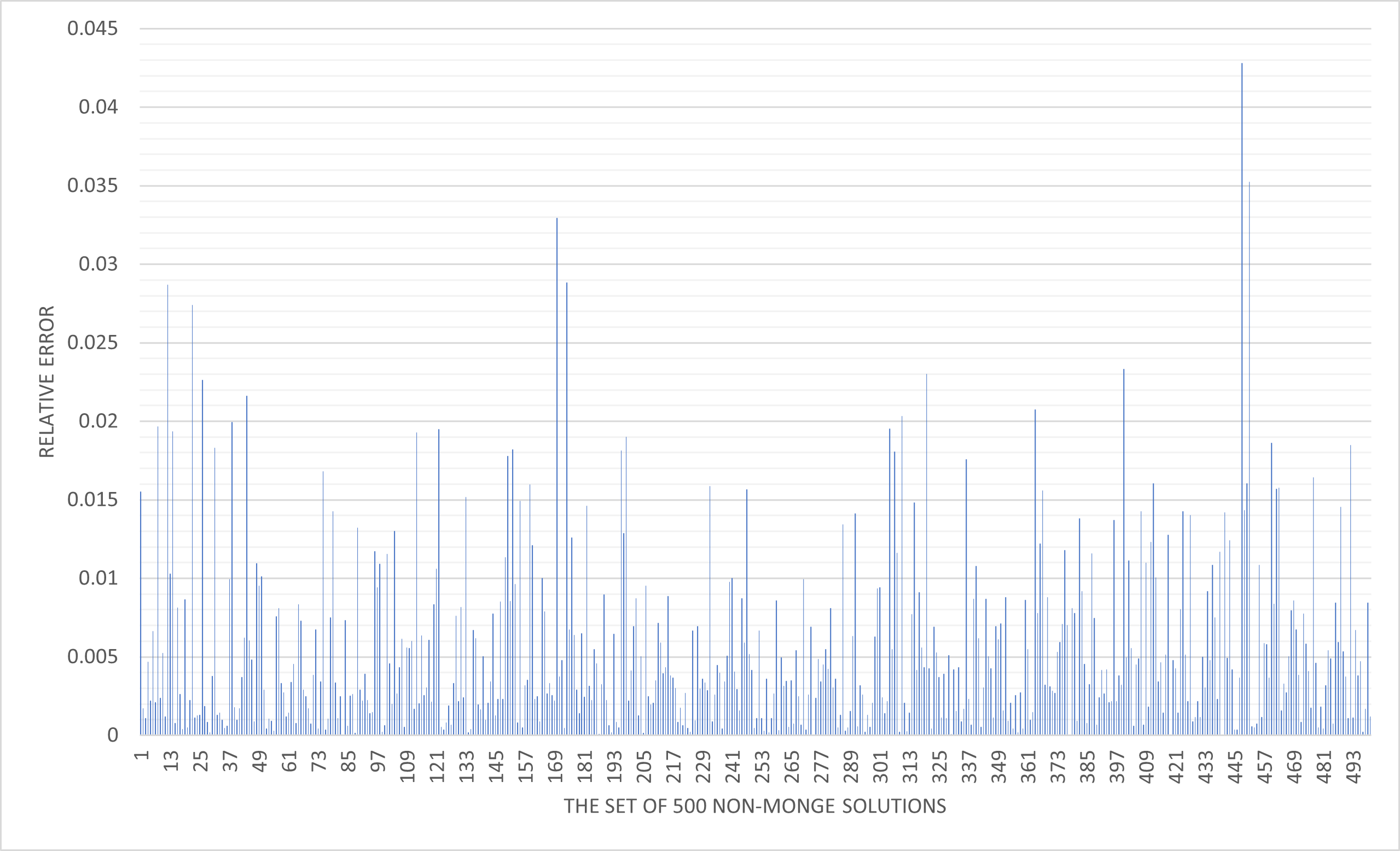}
\caption{The distribution of the errors of minimal Monge costs for $500$ occurrences of multi-marginal
optimal transport problems with $3$-empirical marginals in $\mathbb{R}^3$ which do not have Monge
solutions.}
\label{fig: frequency of non-Monge solutions}
\end{figure}

The counterexample in Example \ref{ex: counterexample} was found by brute force, through randomly
generating many collections of empirical marginals and computing the solutions to \eqref{eqn: OT
problem}. A couple of observations can be gleaned from the data obtained by this procedure:
\begin{enumerate}
\item The occurrence of problems that do not have Monge solutions is quite rare; for three measures
supported on three points in $\R^2$, we found that only  about one in every $5000$ to $8000$ problems
does not have a Monge solution.  The frequency of these counterexamples seems to increase with the
underlying dimension; we found that the failure of the Monge ansatz occurs more often in  $\R^3$.
\item Whenever the Monge ansatz fails, the difference between the minimal Monge cost and the actual
optimal cost is quite small. We identified around $500$ cases when the Monge ansatz fails in $\R^2$;
among these the greatest difference between the cost of the minimal Monge cost and the actual optimal
cost was $4.3\%$.  In fact, in most cases the error was much lower than the above upper bound as is shown
in Figure \ref{fig: frequency of non-Monge solutions}.
\end{enumerate}
The last observation suggest rigorously quantifying the error made by the Monge ansatz as an intriguing
line of future research.  It seems possible that the errors are small enough that the Monge ansatz as an
approximation of \eqref{eqn: OT problem} may still be a  useful tool in some applications.

\bibliographystyle{amsplain}
\bibliography{refs}
\end{document}